\documentclass[9pt]{article}
\usepackage[latin1]{inputenc}
\usepackage{amscd}
\usepackage{amsfonts}
\usepackage{amsgen}
\usepackage{amsmath}
\usepackage{amssymb}
\usepackage{amstext}
\usepackage{amsthm}
\usepackage{graphicx}
\usepackage{indentfirst}
\usepackage{latexsym}
\usepackage{makeidx}
\usepackage{mathrsfs}
\usepackage{epsfig}
\usepackage{wrapfig}
\usepackage[all,knot,arc,import,poly]{xy}
\usepackage[usenames]{color}

\newtheorem{Df}{Definition}[section]
\newtheorem{Teo}[Df]{Theorem}

\newtheorem{Lem}[Df]{Lemma}
\newtheorem{Ex}[Df]{Example}
\newtheorem{Obs}[Df]{Remark}

\newtheorem{Que}[Df]{Question}

\newcommand{\n}{\noindent}
\newcommand{\Dem}{\n{\bf Proof:\;\;}}

\newcommand{\bc}{\begin{center}}
\newcommand{\ec}{\end{center}}

\newcommand{\N}{\mathbb{N}}

\newcommand{\vcinco}{\vspace*{0.5cm}}
\newcommand{\vtres}{\vspace*{0.3cm}}

\def \beq { \begin{equation} }
\def \eeq { \end{equation} }

\newcommand{\baf}{\begin{Afi}\nf{\sl . }}
\newcommand{\eaf}{\end{Afi}}

\newcommand{\bdf}{\begin{Df}\nf{\bf .}}
\newcommand{\edf}{\end{Df}}
\newcommand{\bte}{\begin{Th}\nf{\bf .}}
\newcommand{\ete}{\end{Th}}  
\newcommand{\bco}{\begin{Co}\nf{\bf .}}
\newcommand{\eco}{\end{Co}}
\newcommand{\ble}{\begin{Le}\nf{\bf .}}
\newcommand{\ele}{\end{Le}}
\newcommand{\bpr}{\begin{Pro}\nf{\bf .}}
\newcommand{\epr}{\end{Pro}}
\newcommand{\bex}{\begin{Exa}\nf{\bf .} \rm}
\newcommand{\eex}{\end{Exa}}
\newcommand{\bre}{\begin{Rem}\nf{\bf .} \rm}
\newcommand{\ere}{\end{Rem}}
\newcommand{\bfa}{\begin{Fa}\nf{\bf .} \sl}
\newcommand{\efa}{\end{Fa}}
\newcommand{\bqt}{\begin{Que}\nf{\bf .}}
\newcommand{\eqt}{\end{Que}}

\newcommand{\bct}{\begin{Ct}\nf{\bf .} \rm}
\newcommand{\ect}{\end{Ct}}
\newcommand{\bxa}{\begin{Exa}\nf{\bf .} \rm}
\newcommand{\exa}{\end{Exa}}

\newcommand{\sub}{\subseteq}

\newcommand{\sig}{\mathbb{S}ig}
\newcommand{\cat}{\mathbb{C}at}
\newcommand{\set}{\mathbb{S}et}

\newcommand{\inst}{{\bf Inst}}

\begin{document}
\title{Remarks on Propositional Logics and the categorial relationship between Institutions and $\Pi$-Institutions}
\author{Darllan Concei\c c\~ao\ Pinto \thanks{Institute of Mathematics, Federal University of Bahia, Brazil, darllan\underline{ }math@hotmail.com}, Hugo Luiz Mariano \thanks{Institute of Mathematics and Statistics, University of S\~ao Paulo, Brazil, hugomar@ime.usp.br}}

\maketitle

\begin{abstract}
In this work we explore some applications of the notions of \emph{Institution} and $\pi$-\emph{Institution} in the setting of {\em propositional logics} and establish a precise categorial relation between these notions, i.e., we provide a pair of functors that establishes an adjunction between the categories \inst\  and $\pi$-\inst.
\end{abstract}

\section{Introduction}

The notion of \emph{Institution} was introduced for the first time by Goguen and Burstall in \cite{GB}. This concept formalizes the notion of logical system into a mathematical object, i.e., it provides a \emph{``...categorical abstract model theory which formalizes the intuitive notion of logical system, including syntax, semantic, and satisfaction relation between them...''} \cite{Diac}. This means that it encompasses the abstract concept of universal model theory for a logic. The main (model-theoretical) characteristic is that an institution contains a satisfaction relation between models and sentences that are coherent  under change of notation. First-order (infinitary) logics with Tarski's semantics are natural examples of institutions (see section 2.1).

A variation of the formalism of institutions, the notion of  $\pi$-\emph{Institution}, were defined by Fiadeiro and Sernadas in  \cite{FS} providing an alternative (proof-theoretical)  approach to deductive system \emph{``...replace the notion of model and satisfaction by a primitive consequence operator (\`a la Tarski)''} \cite{FS}. Natural categories of propositional logics (see section 2.2) provide examples of $\pi$-institutions.

In  \cite{FS} and \cite{Vou}   was established a relation between institutions and $\pi$-institutions. On the best of our knowledge, there is no literature on categorial connections between the category of institutions and the category of $\pi$-institutions. In the section 2.3 of the present work,  we provide a precise categorial relationship between these notions, that extends the above mentioned  relation between objects of those categories, more precisely, we determine a pair of adjoint functors between those categories. We finish this work with some remarks concerning, mainly, applications of these tools to the propositional logic setting.

\section{The categories $\inst$ and $\pi-\inst$}

We start giving the definition of institution and $\pi$-institution with their respective notions of morphisms (and comorphisms), and consequently their categories.

\subsection{Institution and its category}
\begin{Df}
An Institution $I=(\mathbb{S}ig,Sen,Mod,\models)$ consists of

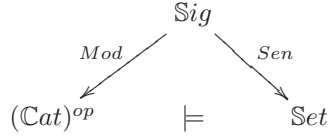
\[\xymatrix{
&\mathbb{S}ig\ar[ld]_{Mod}\ar[rd]^{Sen}&\\
(\mathbb{C}at)^{op}&\models&\mathbb{S}et
}\]

\begin{enumerate}
\item[1.]a category $\mathbb{S}ig$, whose the objects are called {\em signature},
\item[2.]a functor $Sen:\sig\to\set$, for each signature a set whose elements are called {\em sentence} over the signature
\item[3.]a functor $Mod:(\sig)^{op}\to\cat$, for each signature a category whose the objects are called {\em model},
\item[4.]a relation $\models_{\Sigma}\subseteq|Mod(\Sigma)|\times Sen(\Sigma)$ for each $\Sigma\in|\sig|$, called $\Sigma$-{\em satisfaction}, such that for each morphism $h:\Sigma\to\Sigma'$, the compatibility  condition
    \[M'\models_{\Sigma'}Sen(h)(\phi)\ if\ and\ only\ if\ Mod(h)(M')\models_{\Sigma}\phi\]
  holds for each $M'\in|Mod(\Sigma')|$ and $\phi\in Sen(\Sigma)$
\end{enumerate}
\end{Df}

\begin{Ex} \label{inst-examples}
Let $Lang$ denote the  category of  languages $L = ((F_n)_{n \in \N},(R_n)_{n \in \N})$, -- where  $F_n$ is a set of symbols of $n$-ary function symbols and $R_n$ is a   set of symbols of $n$-ary relation symbols, $n \geq 0$ -- and language morphisms\footnote{That can be chosen ``strict" (i.e., $F_n\mapsto F_n'$, $R_n \mapsto R'_n$) or chosen be ``flexible" (i.e., $F_n\mapsto \{n-ary-terms(L')\}$, $R_n \mapsto \{n-ary-atomic-formulas(L')\}$).}. For each pair of cardinals $\aleph_0 \leq \kappa, \lambda \leq \infty$, the category $Lang$  endowed with the usual notion of  $L_{\kappa,\lambda}$-sentences (=  $L_{\kappa,\lambda}$-formulas with no free variable), with the usual association of category of structures and  with the usual (tarskian) notion of satisfaction, gives rise to an institution $I({\kappa,\lambda})$.

\end{Ex}

\begin{Df} Let $I$ and $I'$ be institutions.

(a) An Institution {\bf morphism} $h = (\Phi,\alpha,\beta) : I \to I'$ consists of

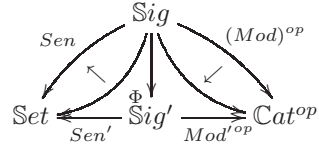
\[\xymatrix{
&\mathbb{S}ig\ar@/^1pc/[ld]_{\nwarrow}\ar@/_/[ld]_{Sen}\ar@/^/[rd]^{(Mod)^{op}}\ar@/_1pc/^{\swarrow}[rd]\ar[d]_>{\Phi}&\\
\set&\sig'\ar[l]^{Sen'}\ar[r]_{{Mod'}^{op}}&\cat^{op}
}\]

\begin{enumerate}
\item[1.]a functor $\Phi:\sig\to\sig'$
\item[2.]a natural transformation $\alpha:Sen'\circ \Phi\Rightarrow Sen$
\item[3.]a natural transformation $\beta:Mod\Rightarrow Mod'\circ\Phi^{op}$
\end{enumerate}
such that the following {\em compatibility condition} holds:
\[m\models_{\Sigma}\alpha_{\Sigma}(\varphi')\ \ if\!f\ \ \beta_{\Sigma}(m)\models'_{\Phi(\Sigma)}\varphi'\]

For any $\Sigma\in\sig$, any $\Sigma$-model $m$ and any $\Phi(\Sigma)$-sentence $\varphi'$.
\vcinco

(b) A triple $f=\langle \phi,\alpha,\beta\rangle:I\to I'$ is a {\bf comorphism} between the given institutions if the following conditions hold:
\begin{itemize}
\item[$\bullet$]$\phi:\sig\to\sig'$ is a functor.
\item[$\bullet$]$\alpha:Sen\Rightarrow Sen'\circ\phi$ and $\beta:Mod'\circ\phi^{op}\Rightarrow Mod$ are natural transformations such that satisfy:

\[m'\models'_{\phi(\Sigma)}\alpha_{\Sigma}(\varphi)\ iff\ \beta_{\Sigma}(m')\models_{\Sigma}\varphi\]

For any $\Sigma\in\sig$,  $m'\in Mod'(\phi(\Sigma))$ and $\varphi\in Sen(\Sigma)$.
\end{itemize}

\end{Df}

\begin{Ex} \label{inst-morph-examples} Given two pairs of cardinals  $(\kappa_i, \lambda_i)$, with  $\aleph_0 \leq \kappa_i, \lambda_i \leq \infty$, $i =0,1$, such that $\kappa_0\leq \kappa_1$ and $\lambda_0  \leq \lambda_1$, then it is induced a morphism and a comorphism of institutions
$(\Phi, \alpha, \beta): I(\kappa_0, \lambda_0) \to I(\kappa_1, \lambda_1)$, given by the {\em same data}: $\sig_{0} = Lang = \sig_{1}$, $Mod_0 = Mod_1 : (Lang)^{op} \to \cat$, $Sen_{i}=L_{\kappa_{i},\lambda_{i}},\ i=0,1$, $\Phi = Id_{Lang} : \sig_0 \to \sig_1$, $\beta := Id : Mod_{i} \Rightarrow Mod_{1-i}$,   $\alpha := inclusion : Sen_{0} \Rightarrow Sen_{1}$.

\end{Ex}

Given $f:I\to I'$ and $f':I\to I''$ comorphisms of institutions, then $f'\bullet f := \langle\phi'\circ\phi, \alpha'\bullet\alpha,\beta'\bullet\beta\rangle$ defines a comorphism $f'\bullet f : I \to I''$, where $(\alpha'\bullet\alpha)_{\Sigma}=\alpha'_{\phi(\Sigma)}\circ\alpha_{\Sigma}$ and $(\beta'\bullet\beta)_{\Sigma}=\beta_{\Sigma}\circ\beta'_{\phi(\Sigma)}$. Let $Id_{I} := \langle Id_{\sig},Id,Id\rangle:I\to I$. It is straitforward to check that these data determines a category\footnote{As usual in category theory, the set theoretical size issues on such global constructions of categories can be addressed by the use of, at least, two Grothendieck's universes.}. We will denote by ${\bf Inst}$ this  category of institutions where the arrows are {\bf comorphisms} of institutions. Of course, it can also be formed a category whose objects are institutions and the arrows are {\bf morphisms} of institutions, but that will be less important here.

\subsection{$\pi$-Institution and its category}

\begin{Df}
A $\pi$-Institution $J=\langle\sig,Sen,\{C_{\Sigma}\}_{\Sigma\in|\sig|}\rangle$ is a triple with its first two components exactly the same as the first two components of an institution and, for every $\Sigma\in|\sig|$, a closure operator $C_{\Sigma}:\mathcal{P}(Sen(\Sigma))\to\mathcal{P}(Sen(\Sigma))$, such that the following coherence conditions holds, for every $f:\Sigma_{1}\to\Sigma_{2}\in Mor(\sig)$:

\[Sen(f)(C_{\Sigma_{1}}(\Gamma))\subseteq C_{\Sigma_{2}}(Sen(f)(\Gamma)),\ for\ all\ \Gamma\subseteq Sen(\Sigma_{1}).\]

\end{Df}

\begin{Df}

Let $J$ and $J'$ be $\pi$-institutions, $g=\langle\phi,\alpha\rangle:J\to J'$ is a \textbf{comorphism} between $\pi$-institution when the following conditions hold:

\begin{itemize}
\item[$\bullet$]$\phi:\sig\to\sig'$ is a functor
\item[$\bullet$]$\alpha:Sen\Rightarrow Sen'\circ\phi$ is a natural transformation such that satisfies the compatibility condition:

\[\varphi\in C_{\Sigma}(\Gamma)\Rightarrow \alpha_{\Sigma}(\varphi)\in C_{\phi(\Sigma)}(\alpha_{\Sigma}(\Gamma))\ for\ all\ \Gamma\cup\{\varphi\}\subseteq\sig(\Sigma)\]
\end{itemize}

\end{Df}

Let $g:J\to J'$ and $g':J'\to J''$ be comorphisms of $\pi$-institutions. $g'\bullet g$ is defined as the two first components of composition of comorphisms of institutions. The identity (co)morphism is given as the two first components of the comorphism identity of institution. We will denote by $\pi$-{\bf Inst} the category of $\pi$-institutions and with arrows its comorphisms.
\vtres

\begin{Ex} \label{pi-inst-examples} In \cite{AFLM}, \cite{FC} and \cite{MaMe} are considered some categories of propositional logics -- $l = (\Sigma, \vdash)$,  where $\Sigma = (\Sigma_n)_{n \in \N}$  is a finitary signature and $\vdash \subseteq P(Form(\Sigma))\times Form(\Sigma))$ is a tarskian consequence operator-- with morphisms, $f : (\Sigma, \vdash) \rightarrow (\Sigma', \vdash')$, some kind of signature morphism $f : \Sigma \rightarrow \Sigma'$ --``strict"  or ``flexible"-- that induces a translation or interpretation : $\Gamma \cup \{\psi\} \subseteq Form(\Sigma)$, $\Gamma \vdash {\psi} \Rightarrow \check{f}[\Gamma] \vdash'  \check{f}(\psi)$.

(a)  To the category of propositional logics endowed with ``{\em flexible} morphisms" ${\cal L}_f$ (respectively, endowed with ``{\em strict} morphisms" ${\cal L}_s$) is associated an $\pi$-institution  $J_f$ (respectively, $J_s$) in the following way:\\
$\bullet$ $\sig_f$ := ${\cal L}_f$;\\
$\bullet$  $Sen_f : \sig_f \to \set$, given by $(f: (\Sigma, \vdash) \to (\Sigma', \vdash)) \ \mapsto\ (\check{f} : F_\Sigma(X) \to F_{\Sigma'}(X))$;\\
$\bullet$ For each $l = (\Sigma, \vdash) \in |\sig_f|$, $C_{l} : P(F_\Sigma(X)) \to P(F_\Sigma(X))$ is given by $C_{l}(\Gamma):= \{\phi \in  F_\Sigma(X) :\Gamma\vdash_l \phi\}$, for each $\Gamma \sub  F_\Sigma(X)$.

(b) In \cite{MaMe}, is The ``inclusion" functor $(+)_L : {\cal L}_s \to {\cal L}_f$, induces a comorphism (and also a morphism!) of the associated $\pi$-institutions $(+) := ((+)_L, \alpha^+) : J_s \to J_f$, where, for each $l = (\Sigma, \vdash) \in \sig_s = {\cal L}_s$, $\alpha^+(l) = Id_{F_{\Sigma}(X)} : F_{\Sigma}(X) \to F_{\Sigma}(X)$.

\end{Ex}

\subsection{An adjunction between \inst\ and $\pi$-\inst}

In order to establish the adjunction between $Inst$ and $\pi-Inst$ we introduce the following:

Let $I=\langle\sig,Sen,Mod,\models\rangle$ be an institution. Given $\Sigma\in|\sig|$, consider
\[\Gamma^{\star}=\{m\in Mod(\Sigma);\ m\models_{\Sigma} \varphi\ for\ all\ \varphi \in\Gamma\}\ \ and\]
\[M^{\star}=\{\varphi\in Sen(\Sigma);\ m\models_{\Sigma}\varphi\ for\ all\ m\in M\}\]

for any $\Gamma\subseteq Sen(\Sigma)$ and $M\subseteq Mod(\Sigma)$. Clearly, these mappings establishes a Galois connection. Thus $C^{I}_{\Sigma}(\Gamma) : =\Gamma^{\star\star}$, defines a closure operator for any $\Sigma\in|\sig|$ (\cite{Vou}).

 The following lemma describes the behavior of these Galois connections through institutions comorphisms.

\begin{Lem}\label{lema1}
Let $f=\langle\phi,\alpha,\beta\rangle:I\to I'$ an arrow in {\bf Inst}. Then given $\Gamma\subseteq Sen(\Sigma)$ and $M\subseteq|Mod(\phi(\Sigma))|$, the following conditions holds:
\begin{itemize}
\item[1)]$\beta_{\Sigma}[(\alpha_\Sigma[\Gamma])^{\star}]\subseteq \Gamma^{\star}$
\item[2)]$\alpha_{\Sigma}[(\beta_{\Sigma}[M])^{\star}]\subseteq M^{\star}$
\end{itemize}
\end{Lem}

\Dem \emph{1)} Let $m\in\beta_{\Sigma}[(\alpha_\Sigma[\Gamma])^{\star}]$. So there is $m'\in\alpha_\Sigma[{\Gamma}]^{\star}$ such that $\beta_{\Sigma}(m')=m$. As $m'\in\alpha_{\Sigma}[\Gamma]^{\star}$, hence $m'\models'_{\phi(\Sigma)}\alpha_{\Sigma}[\Gamma]\ \Leftrightarrow\ \beta_{\Sigma}(m')\models_{\Sigma}\Gamma\ \Leftrightarrow\ m\models_{\Sigma}\Gamma$. Then $m\in\Gamma^{\star}$.

\emph{2)} Let $\varphi\in\alpha_{\Sigma}[(\beta_{\Sigma}[M])^{\star}]$. So there is $\psi\in\beta_{\Sigma}[M]^{\star}$such that $\alpha_{\Sigma}(\psi)=\varphi$. Since $\psi\in(\beta_{\Sigma}[M])^{\star}$, hence $\beta_{\Sigma}[m]\models_{\Sigma}\psi\ \Leftrightarrow\ m\models_{\phi(\Sigma)}\alpha_{\Sigma}(\psi)\ \Leftrightarrow\ m\models_{\phi(\Sigma)}\varphi$ for any $m\in M$. Therefore $\varphi\in M^{\star}$.
\qed
\vcinco

Define the following application:

$\begin{array}{rrcl}
F:&\inst&\longrightarrow&\pi-\inst\\
&I&\longmapsto&F(I)=\langle\sig,Sen,\{C^{I}_{\Sigma}\}_{\Sigma\in|\sig|}\rangle
\end{array}
$

In order to provide the well-definition of $F$, it is enough to prove the compatibility condition for $\{C^{I}_{\Sigma}\}_{\Sigma\in|\sig|}$, i.e., given $f:\Sigma_1\to\Sigma_2$ and $\Gamma\subseteq Sen(\Sigma_1)$, then $Sen(f)(C^{I}_{\Sigma_1}(\Gamma))\subseteq C^{I}_{\Sigma_2}(Sen(f)(\Gamma))$. Let $\varphi_2\in Sen(f)(C^{I}_{\Sigma_1}(\Gamma))$, then there is $\varphi_1\in\Gamma^{**}$ such that $Sen(f)(\varphi_1)=\varphi_2$. Let $m\in(Sen(f)(\Gamma))^{*}$. So $m\models_{\Sigma_2}Sen(f)(\Gamma)$. By compatibility condition in institutions we have that $Mod(f)(m)\models_{\Sigma_1}\Gamma$, thus $Mod(f)(m)\in\Gamma^{*}$. Since $\varphi_1\in\Gamma^{**}$ we have that $Mod(f)(m)\models_{\Sigma_1}\varphi_1$, hence $m\models_{\Sigma_2}Sen(f)(\varphi_1)=\varphi_2$. Therefore $\varphi_2\in (Sen(f)(\Gamma))^{**}=C^{I}_{\Sigma_2}(Sen(f)(\Gamma))$.

Now let $f=\langle\phi,\alpha,\beta\rangle:I\to I'$ be a comorphism of institutions. Then consider $F(f)=\langle\phi,\alpha\rangle$. Notice that $F(f)$ is a comorphism between $F(I)$ and $F(I')$. Indeed, it is enough to prove that $F(f)$ satisfies the compatibility condition. Let $\Gamma\cup\{\varphi\}\subseteq Sen(\Sigma)$ for some $\Sigma\in|\sig|$. Suppose that $\alpha_{\Sigma}(\varphi)\not\in C^{I}_{\phi(\Sigma)}(\alpha_{\Sigma}[\Gamma])$. Hence $\alpha_{\Sigma}(\varphi)\not\in\alpha_{\Sigma}[\Gamma]^{\star\star}$. Therefore $\alpha_{\Sigma}[\Gamma]^{\star}\not\models'_{\phi(\Sigma)}\alpha_{\Sigma}(\alpha)$. Thus there is $m\in \alpha_{\Sigma}[\Gamma]^{\star}$ such that $m\not\models'_{\phi(\Sigma)}\alpha_{\Sigma}(\varphi)$. Hence $\beta_{\Sigma}(m)\not\models_{\Sigma}\varphi$. Due to \ref{lema1} 1) we have that $\beta_{\Sigma}(m)\in\Gamma^{\star}$. Therefore $\varphi\not\in\Gamma^{\star\star}=C^{I}_{\Sigma}(\Gamma)$.

Now let $f:I\to I'$ and $f':I'\to I''$ comorphism of institutions. $F(f'\bullet f)=\langle\phi'\circ\phi,\alpha'\bullet\alpha\rangle=F(f')\bullet F(f)$ and $F(Id_{I})=Id_{F(I)}$. Then $F$ is a functor.
\vtres

Consider now the application:

$\begin{array}{rrcl}
G:&\pi-\inst&\longrightarrow&\inst\\
&J&\longmapsto&G(J)=\langle\sig,Sen,Mod^{J},\models^{J}\rangle
\end{array}
$

Where:\\

$\bullet$ The two first components of the $\pi-$institution are preserved.

$\bullet$ $Mod^{J}:\sig\to\cat^{op}$.\\
 $Mod^{J}(\Sigma) : =\{C_{\Sigma}(\Gamma);\ \Gamma\subseteq Sen(\Sigma)\} \subseteq P(Sen(\Sigma))$ is viewed as a poset category and, given $f:\Sigma\to\Sigma'$, $Mod^{J}(f)=Sen(f)^{-1}$. \\
$Mod^{J}(f)$ is well defined. Indeed: Let $\Gamma\subseteq Sen(\Sigma')$ and $\varphi\in C_{\Sigma}(Sen(f)^{-1}(C_{\Sigma'}(\Gamma)))$.
\[\begin{array}{rcl}
Sen(f)(\varphi)\in Sen(f)[C_{\Sigma}(Sen(f)^{-1}(C_{\Sigma'}[\Gamma]))]&\subseteq& C_{\Sigma}[Sen(f)(Sen(f)^{-1}(C_{\Sigma}[\Gamma]))]\\
&\subseteq& C_{\Sigma'}(C_{\Sigma'}[\Gamma])=C_{\Sigma'}[\Gamma]
\end{array}\]
Therefore $\varphi\in Sen(f)^{-1}(C_{\Sigma}[\Gamma])$. It is easy to see that $Mod^{J}$ is a contravariant functor.

$\bullet$ Define $\models^{J}\subseteq|Mod(\Sigma)|\times Sen(\Sigma)$ as a relation such that given $m\in Mod(\Sigma)$ and $\varphi\in Sen(\Sigma)$, $m\models^{J}_{\Sigma}\varphi$ if and only if $\varphi\in m$. Let $f:\Sigma\to\Sigma'$, $\varphi\in Sen(\Sigma)$ and $m'\in|Mod(\Sigma')|$.

$\begin{array}{rcl}
Mod^{J}(f)(m')\models^{J}_{\Sigma}\varphi&\Leftrightarrow&Sen(f)^{-1}(m')\models^{J}_{\Sigma}\varphi\\
&\Leftrightarrow&\varphi\in Sen(f)^{-1}(m')\\
&\Leftrightarrow& Sen(f)(\varphi)\in m'\\
&\Leftrightarrow& m'\models^{J}_{\Sigma'}Sen(f)(\varphi)
\end{array}$

Therefore the compatibility condition is satisfied and then we have that $G(J)$ is an institution.
\\

Now let $h=\langle\phi,\alpha\rangle:J\to J'$ be a comorphism of $\pi$-institution. Define for any $\Sigma\in|\sig|$ $\beta_{\Sigma}:Mod^{J'}\circ\phi(\Sigma)\to Mod^{J}(\Sigma)$ where $\beta_{\Sigma}(m)=\alpha_{\Sigma}^{-1}(m)$. We prove that $\beta_{\Sigma}$ is well defined, i.e., $\alpha_{\Sigma}^{-1}(m)\in Mod^{J}(\Sigma)$. Let $\varphi\in C_{\Sigma}(\alpha_{\Sigma}^{-1}(m))$. Since $h$ is a morphism of $\pi$-institution, then $\alpha_{\Sigma}(\varphi)\in C_{\phi(\Sigma)}(\alpha_{\Sigma}(\alpha_{\Sigma}^{-1}(m)))\subseteq C_{\phi(\Sigma)}(m)=m$. Therefore $\varphi\in\alpha_{\Sigma}^{-1}(m)$.

Now we prove that $\beta$ is a natural transformation. Let $f:\Sigma_{1}\to\Sigma_{2}$. Since $\alpha$ is a natural transformation, the following diagram commutes:

\[
\xymatrix{
P(Sen(\Sigma_{1}))&P(Sen'(\phi(\Sigma_{1})))\ar[l]_{\alpha_{\Sigma_{1}}^{-1}}\\
P(Sen(\Sigma_{2}))\ar[u]^{Sen(f)^{-1}}&P(Sen'(\phi(\Sigma_{2})))\ar[u]_{Sen'(\phi(f))^{-1}}\ar[l]^{\alpha_{\Sigma_{2}}^{-1}}
}\]

Using this commutative diagram we are able to prove that the following diagram commutes:

\[
\xymatrix{
Mod^{J'}\circ\phi(\Sigma_{1})\ar[r]^{\beta_{\Sigma_{1}}}&Mod^{J}(\Sigma_{1})\\
Mod^{J'}\circ\phi(\Sigma_{2})\ar[u]^{Mod^{J'}(\phi(f))}\ar[r]_{\beta_{\Sigma_{2}}}&Mod^{J}(\Sigma_{2})\ar[u]_{Mod^{J}(f)}
}\]

Let $m\in Mod^{J'}\circ\phi(\Sigma_{2})$.

$
\begin{array}{rcl}
Mod^{J}(f)\circ\beta_{\Sigma_{2}}(m)&=&Mod^{J}(f)(\alpha^{-1}_{\Sigma_{2}}(m))\\
&=&Sen(f)^{-1}(\alpha^{-1}_{\Sigma_{2}}(m))\\
&=&\alpha_{\Sigma_{1}}^{-1}(Sen(\phi(f))^{-1}(m))\\
&=&\beta_{\Sigma_{1}}(Sen(\phi(f))^{-1}(m))\\
&=&\beta_{\Sigma_{1}}\circ Mod^{J'}(\phi(f))(m)
\end{array}
$
\\
$G(h)=\langle\phi,\alpha,\beta\rangle$ is a comorphism of institution. Indeed, it is enough to prove the compatibility condition. Let $m\in Mod^{J'}(\phi(\Sigma))$ and $\varphi\in Sen(\Sigma)$.

$
\begin{array}{rcl}
m\models^{J'}_{\phi(\Sigma)}\varphi \alpha_{\Sigma}(\varphi)&\Leftrightarrow&\alpha_{\Sigma}(\varphi)\in m\\
&\Leftrightarrow&\varphi\in\alpha_{\Sigma}^{-1}(m)\\
&\Leftrightarrow&\varphi\in\beta_{\Sigma}(m)\\
&\Leftrightarrow&\beta_{\Sigma}(m)\models^{J}_{\Sigma}(m)\varphi
\end{array}$

It is easy to see that $G$ is a functor.

\begin{Teo}
The functors $F : \inst \rightarrow \pi-\inst$ and $G : \pi-\inst \rightarrow \inst$ defined above establish an adjunction $G \dashv F$ between the categories $\inst$ and $\pi-\inst$.
\end{Teo}

\Dem

Define the application $\eta_{J}=\langle Id_{\sig},Id_{Sen}\rangle:J\to F(G(J))$ for each $\pi$-Institution $J=\langle\sig,Sen,\{C_{\Sigma}\}_{\Sigma\in|\sig|}\rangle$. This application is well defined. Indeed, we prove that $C_{\Sigma}=C^{G(I)}_{\Sigma}$ for any $\Sigma\in|\sig|$. By definition of the functor $G$, notice that given $\Sigma\in|\sig|$ and $\Gamma\subseteq Sen(\Sigma)$, $C_{\Sigma}(\Gamma)\in \Gamma^{\star}=\{m\in Mod(\Sigma);\ m\models^{J}_{\Sigma}\Gamma\}$. Moreover $C_{\Sigma}(\Gamma)\subseteq m$ for every $m\in \Gamma^{\star}$. Then for any $\varphi\in Sen(\Sigma)$

$\begin{array}{rcl}
\varphi\in C_{\Sigma}(\Gamma)&\Leftrightarrow&\varphi\in m\ for\ all\ m\in\Gamma^{\star}\\
&\Leftrightarrow&m\models^{J}_{\Sigma}\varphi\ for\ all\ m\in\Gamma^{\star}\\
&\Leftrightarrow&\varphi\in\Gamma^{\star\star}=\{\psi\in Sen(\Sigma);\ \Gamma^{\star}\models^{J}_{\Sigma}\psi\}\\
&\Leftrightarrow&\varphi\in C^{G(J)}_{\Sigma}(\Gamma).
\end{array}$

It is clear that $(\eta_{J})_{J\in|\pi-\inst|}$ is a natural transformation. It remains to prove that $\eta_{J}$ satisfies the universal property for any $J\in|\pi-\inst|$.

Let $h=\langle\phi,\alpha\rangle:J\to F(I)$ where $J=\langle\sig,Sen,\{C_{\Sigma}\}_{\Sigma\in|\sig|}\rangle$ is a $\pi-$institution, $I=\langle\sig',Sen',Mod',\models'\rangle$ an institution and $h$ a morphism of $\pi-$institution. Define $\bar{h}=\langle\phi,\alpha,\beta\rangle:G(J)\to I$ where the first two components are the same of $h$ and given $\Sigma\in|\sig|$, $\beta_{\Sigma}:Mod'\circ\phi(\Sigma)\to Mod^{J}(\Sigma)$ such that $\beta_{\Sigma}(m)=\alpha_{\Sigma}^{-1}[m^{\star}]$. $\beta_{\Sigma}$ is well defined. Indeed, notice that $m^{\star}=m^{\star\star\star}$ for any $m\in Mod'(\phi(\Sigma))$. Since $C^{I}_{\Sigma}(\Gamma)=\Gamma^{\star\star}$, therefore $m^{\star}=C^{I}_{\Sigma}(m^{\star})$. We have shown that as $h$ is a morphism of $\pi-$institution, $\alpha^{-1}_{\Sigma}(m^{\star})=\alpha^{-1}_{\Sigma}(C^{I}_{\Sigma}(m^{\star}))\in Mod^{J}$.

Now we prove that $(\beta_{\Sigma})_{\Sigma\in|\sig|}$ is a natural transformation. Let $f:\Sigma_{1}\to \Sigma_{2}$. Then given $m\in Mod'\circ\phi(\Sigma_{2})$

\[\xymatrix{
Mod'\circ\phi(\Sigma_{1})\ar[r]^{\beta_{\Sigma_{1}}}&Mod^{J}(\Sigma_{1})\\
Mod'\circ\phi(\Sigma_{2})\ar[r]_{\beta_{\Sigma_{2}}}\ar[u]^{Mod'(\phi(f))}&Mod^{J}(\Sigma_{2})\ar[u]_{Mod^{J}(f)}
}\]

$\begin{array}{rcl}
Mod^{J}(f)(\beta_{\Sigma_{2}}(m))&=&Sen(f)^{-1}(\alpha^{-1}_{\Sigma_{2}}(m^{\star}))\\
&=&\alpha^{-1}_{\Sigma_{1}}(Sen(\phi(f)^{-1})(m^{\star}))\\
&=_\dagger &\alpha^{-1}_{\Sigma_{1}}((Mod(\phi(f))(m^{\star}))^{\star})\\
&=&\beta_{\Sigma_{1}}(Mod(\phi(f))(m)).
\end{array}$

The justification of the equality $(\dagger)$ is:

$\begin{array}{rcl}
\varphi\in Sen(\phi(f))^{-1}(m^{\star})&\Leftrightarrow& Sen(\phi(f))(\varphi)\in m^{\star}\\
&\Leftrightarrow&m\models_{\phi(\Sigma_{2})}Sen(\phi(f))(\varphi)\\
&\Leftrightarrow&Mod(\phi(f))(m)\models\Sigma_{2}\varphi\\
&\Leftrightarrow&\varphi\in(Mod(\phi(f))(m))^{\star}
\end{array}
$

Hence $\beta$ is a natural transformation. Therefore $\bar{h}$ is a comorphism between $G(I)$ and $I$. Observe that $F(\bar{h})=\langle\phi,\alpha\rangle=h$. Then we have the following diagram commuting:

\[\xymatrix{
J\ar[r]^{\eta_{J}}\ar[dr]_{h}&F(G(J))\ar[d]^{F(\bar{h})}\\
&F(I)
}
\]

Moreover, clearly $\bar{h}$ is the unique arrow such that the diagram above commutes. Hence $G\dashv F$.\qed

\begin{Obs} \label{adj-rem} Note that $F \circ G = Id_{\pi-\inst}$ and the unity of this adjunction, the natural transformation $\eta : Id_{\pi-\inst} \rightarrow F\circ G$, is the identity. Thus the category ${\pi-\inst}$ can be seen as a full co-reflective subcategory of $\inst$.

\end{Obs}









\section{Final remarks and future works}

\begin{Obs} \label{catlog-rem} In \cite{MaMe} it is presented a right adjoint $(-)_L : {\cal L}_f \to {\cal L}_s$ to the ``inclusion" functor $(+)_L : {\cal L}_s \to {\cal L}_f$ (see Example \ref{pi-inst-examples}).  It will be interesting  understand the role of these adjoint pair of functors  between the logical categories (${\cal L}_s, {\cal L}_s$) at the $\pi$-institutional level ($J_f, J_s$).

\end{Obs}

\begin{Obs} \label{Gliv-rem}  The ``proof-theoretical" Example \ref{pi-inst-examples}, that provides  {\em $\pi$-institutions} ($J_f, J_s$) for a categories of propositional logics (${\cal L}_s, {\cal L}_s$), lead us to search  an analogous ``model-theoretical" version of it that is different from the canonical one (i.e., that  obtained by applying the functor $G : \pi-\inst \rightarrow \inst$): In \cite{MaPi2}, we   provide (another) \emph{institutions} for each   category of propositional logics, through the use of the notion of a {\em matrix} for a propositional logic. Moreover, by a convenient modification of this later construction, we provide in \cite{MaPi2} an institution for each ``equivalence class" of algebraizable logic: this  enable us to  apply notions and results from Institution Theory in the propositional logic setting  and derive, from the introduction of the notion of ``Glivenko's context", a strong and general form of Glivenko's Theorem relating two ``well-behaved" logics.
\end{Obs}

The examination of the content mentioned in both the remarks above could lead naturally  to consider  new  categories of propositional logics and to a new notions of morphism of ($\pi$-)institutions.

This work also open a way to investigate categorial properties of the categories of institutions and $\pi$-institutions with many kinds of morphisms in each of them.


\end{document}